\documentclass[a4paper,10pt,reqno]{amsart}
\usepackage[a4paper,tmargin=40mm,bmargin=35mm,hmargin=30mm]{geometry}
\usepackage[T1]{fontenc}
\usepackage{lmodern}
\usepackage{amsmath}
\usepackage{amssymb}
\usepackage{amsthm}
\usepackage{stmaryrd}
\usepackage{tikz}
\usetikzlibrary{cd}
\usepackage{booktabs}
\usepackage{multirow}
\usepackage{multicol}
\usepackage{here}
\usepackage{aliascnt}
\usepackage{hyperref}
\usepackage[noabbrev]{cleveref}

\newcommand{\NewTheorem}[2]{
	\newaliascnt{#1}{TheoremEnvironment}
	\newtheorem{#1}[#1]{#1}
	\aliascntresetthe{#1}
	\crefname{#1}{#1}{#2}
	\Crefname{#1}{#1}{#2}
}

\theoremstyle{definition}

\NewTheorem{Definition}{Definitions}
\NewTheorem{Remark}{Remarks}
\NewTheorem{Axiom}{Axioms}
\NewTheorem{Example}{Examples}
\NewTheorem{Observation}{Observations}
\NewTheorem{Convention}{Conventions}
\NewTheorem{Notation}{Notations}
\NewTheorem{Setting}{Settings}
\NewTheorem{Question}{Questions}
\NewTheorem{Answer}{Answers}
\NewTheorem{Conjecture}{Conjectures}
\NewTheorem{Problem}{Problems}
\NewTheorem{Solution}{Solutions}
\NewTheorem{Goal}{Goals}
\NewTheorem{Comment}{Comments}
\NewTheorem{Aim}{Aims}
\NewTheorem{Caution}{Cautions}
\NewTheorem{Exercise}{Exercises}

\theoremstyle{plain}
\NewTheorem{Proposition}{Propositions}
\NewTheorem{Lemma}{Lemmas}
\NewTheorem{Theorem}{Theorems}
\NewTheorem{Corollary}{Corollaries}

\crefname{enumi}{}{}
\Crefname{enumi}{}{}
\creflabelformat{enumi}{(#2#1#3)}
\crefname{enumii}{}{}
\Crefname{enumii}{}{}
\creflabelformat{enumii}{(#2#1#3)}
\crefname{enumiii}{}{}
\Crefname{enumiii}{}{}
\creflabelformat{enumiii}{(#2#1#3)}

\makeatletter
\renewcommand{\p@enumii}{}
\renewcommand{\p@enumiii}{}
\makeatother

\numberwithin{equation}{section}
\crefname{equation}{}{}
\Crefname{equation}{}{}
\creflabelformat{equation}{(#2#1#3)}

\newcommand{\SwapSymbols}[1]{
	\expandafter\let\expandafter\temporarysymbol\csname #1\endcsname
	\expandafter\let\csname #1\expandafter\endcsname\csname var#1\endcsname
	\expandafter\let\csname var#1\endcsname\temporarysymbol
}

\SwapSymbols{epsilon}
\SwapSymbols{phi}
\SwapSymbols{Gamma}
\SwapSymbols{Delta}
\SwapSymbols{Theta}
\SwapSymbols{Lambda}
\SwapSymbols{Xi}
\SwapSymbols{Pi}
\SwapSymbols{Sigma}
\SwapSymbols{Upsilon}
\SwapSymbols{Phi}
\SwapSymbols{Psi}
\SwapSymbols{Omega}

\newcommand{\bbN}{\mathbb{N}}

\newcommand{\cG}{\mathcal{G}}

\newcommand{\cX}{\mathcal{X}}

\newcommand{\kp}{\mathfrak{p}}

\newcommand{\resp}{resp.\ }

\let\originalleft\left
\let\originalright\right
\renewcommand{\left}{\mathopen{}\mathclose\bgroup\originalleft}
\renewcommand{\right}{\aftergroup\egroup\originalright}

\newcommand{\set}[2][]{\mathopen{#1\{}#2\mathclose{#1\}}}

\newcommand{\setwithtext}[2][]{\mathopen{#1\{}\,\textnormal{#2}\,\mathclose{#1\}}}
\newcommand{\setwithcondition}[3][]{\mathopen{#1\{}\,#2\mathrel{#1|}#3\,\mathclose{#1\}}}

\newcommand{\into}{\hookrightarrow}

\newcommand{\onto}{\twoheadrightarrow}

\newcommand{\isoto}{\xrightarrow{\smash{\raisebox{-0.25em}{$\sim$}}}}

\newcommand{\bc}{\mathbf{c}}
\newcommand{\br}{\mathbf{r}}
\newcommand{\KC}{K\langle\langle C\rangle\rangle}
\newcommand{\KtC}{K\langle\langle\widetilde{C}\rangle\rangle}

\DeclareMathOperator{\id}{id}

\DeclareMathOperator{\Mod}{Mod}

\DeclareMathOperator{\Spec}{Spec}

\DeclareMathOperator{\supp}{supp}

\DeclareMathOperator{\ASpec}{ASpec}

\DeclareMathOperator{\ASupp}{ASupp}

\title[Construction of Grothendieck categories]{Construction of Grothendieck categories with enough compressible objects using colored quivers}
\subjclass[2010]{18E15 (Primary), 16D90, 16P40, 16G30 (Secondary)}
\keywords{Grothendieck category; atom spectrum; Gabriel spectrum; colored quiver; partially ordered set}

\author{Ryo Kanda}
\address[Ryo Kanda]{Department of Mathematics, Graduate School of Science, Osaka University, Toyonaka, Osaka, 560-0043, Japan}
\email{ryo.kanda.math@gmail.com}

\begin{document}

\begin{abstract}
	We introduce a new method to construct a Grothendieck category from a given colored quiver. This is a variant of the construction used to prove that every partially ordered set arises as the atom spectrum of a Grothendieck category. Using the new method, we prove that for every finite partially ordered set, there exists a locally noetherian Grothendieck category such that every nonzero object contains a compressible subobject and its atom spectrum is isomorphic to the given partially ordered set.
\end{abstract}

\maketitle
\tableofcontents

\section{Introduction}
\label{7891873903}

The \emph{atom spectrum} of a Grothendieck category is a generalization of the set of prime ideals of a commutative ring. The idea to associate a space to a Grothendieck category can be found in \cite{MR0232821}, where the \emph{Gabriel spectrum} was defined to be the set of isomorphism classes of indecomposable injective objects. If $R$ is a commutative noetherian ring, Matlis' result ensures that there is a canonical bijection between the Gabriel spectrum of $\Mod R$ and the prime spectrum $\Spec R$, where $\Mod R$ is the category of all $R$-modules. The atom spectrum is a variant of the Gabriel spectrum, which is defined so that there is a canonical bijection between the atom spectrum of $\Mod R$ and $\Spec R$ for an arbitrary commutative ring $R$. The idea of atom spectrum was given by Storrer \cite{MR0360717} to reformulate Goldman's work \cite{MR0245608} on primary decomposition of modules over noncommutative rings. The topological structure on the atom spectrum (or the Gabriel spectrum) has been used to classify localizing subcategories of the given Grothendieck category (see \cite[Proposition~VI.2.4]{MR0232821}, \cite[Theorem~3.8]{MR1434441}, \cite[Corollary~4.3]{MR1426488}, \cite[Theorem~5.5]{MR2964615}, and \cite[Theorem~7.8]{MR3452186}).

In this paper, we will focus on the partial order on the atom spectrum $\ASpec\cG$ of a Grothendieck category $\cG$. The partial order is naturally defined as the \emph{specialization order} with respect to the topology on $\ASpec\cG$, and for a commutative ring $R$, the partial order on $\ASpec(\Mod R)$ is identified with the inclusion of prime ideals via the canonical bijection $\ASpec(\Mod R)\isoto\Spec R$.

The possible partial order structure of the prime spectrum of a commutative ring was completely determined by Hochster \cite[Proposition~10]{MR0251026} and Speed \cite[Corollary~1]{MR0306061}: A partially ordered set is isomorphic to $\Spec R$ for some commutative ring $R$ if and only if it can be written as the inverse limit of finite partially ordered sets in the category of all partially ordered sets. It should be noticed that this property implies several well-known properties of $\Spec R$ such as the fact that every prime ideal is contained in some maximal (prime) ideal. When we only consider commutative noetherian rings, the partial order structure of $\Spec R$ is much more restrictive. It was shown by de Souza Doering and Lequain \cite[Theorem~B]{MR574801} that a \emph{finite} partially ordered set $P$ is isomorphic to $\Spec R$ for some commutative noetherian ring $R$ if and only if there is no chain of the form $x<y<z$ in $P$. Thus such a ring $R$ should have dimension at most one.

The analogous problems for Grothendieck categories have quite different answers. We showed in \cite[Theorem~1.2 (2)]{MR3351569} that \emph{every} partially ordered set arises as the atom spectrum of a Grothendieck category. As an analog to the problem for commutative noetherian rings, we also proved that \emph{every} finite partially ordered set arises as the atom spectrum of some locally noetherian Grothendieck category (\cite[Theorem~1.4 (2)]{MR3351569}). The main idea of these results was to construct a Grothendieck category from a given \emph{colored quiver}.

One of the significant properties of the category $\Mod R$ for a commutative noetherian ring $R$, among all locally noetherian Grothendieck categories, is that it has enough compressible objects. We say that a Grothendieck category $\cG$ has \emph{enough compressible objects} if every nonzero object in $\cG$ contains a compressible subobject (see \cref{0923109234} for the definition of a compressible object). This holds for $\Mod R$ since every nonzero $R$-module contains a submodule isomorphic to $R/\kp$ for some $\kp\in\Spec R$, which is compressible. The category of modules over a right noetherian ring is a locally noetherian Grothendieck category, but it does not necessarily have enough compressible objects (\cref{0981081238}). So it is natural to ask whether having enough compressible objects restricts the partial order structure of the atom spectrum.

The aim of this paper is to show that the answer to this question is \emph{no}. We provide a new method to construct a Grothendieck category from a colored quiver and prove the following result:

\begin{Theorem}[\cref{0893340933}]\label{0984767827}
	Let $P$ be a finite partially ordered set. Then there exists a locally noetherian Grothendieck category $\cG$ satisfying the following properties:
	\begin{enumerate}
		\item\label{0982752493} $\ASpec\cG$ is isomorphic to $P$ as a partially ordered set.
		\item\label{0242394852} $\cG$ has enough compressible objects.
	\end{enumerate}
\end{Theorem}

\subsection*{Acknowledgements}
\label{9172648243}

The author was a JSPS Overseas Research Fellow and supported by JSPS KAKENHI Grant Numbers JP17K14164 and JP16H06337. This work was done while the author was visiting the University of Washington. I would like to thank S.~Paul Smith for his hospitality.

\section{Preliminaries}
\label{0981980374}

In this section, we collect some basic materials on Grothendieck categories and their atom spectra.

\begin{Convention}\label{0930932244}
	A \emph{direct sum} (\resp a \emph{direct product}) in a Grothendieck category means a possibly infinite coproduct (\resp product) whose index set is set-theoretically small.
\end{Convention}

\begin{Definition}\label{0198341923}
	Let $\cG$ be a Grothendieck category.
	\begin{enumerate}
		\item\label{0109812348} A \emph{weakly closed subcategory} (also called a \emph{prelocalizing subcategory}) of $\cG$ is a full subcategory closed under subobjects, quotient, and direct sums.
		\item\label{3781298133} A \emph{closed subcategory} of $\cG$ is a weakly closed subcategory that is also closed under direct products.
	\end{enumerate}
\end{Definition}

\begin{Remark}\label{9023409824}
	For a ring $R$, we denote by $\Mod R$ the category of all (set-theoretically small) right $R$-modules. This is a Grothendieck category. For a two-sided ideal $I$ of $R$, the category $\Mod(R/I)$ is canonically identified with the closed subcategory
	\begin{equation*}
		\setwithcondition{M\in\Mod R}{MI=0}
	\end{equation*}
	of $\Mod R$. It is well known that all closed subcategories of $\Mod R$ is of this form (see, for example, \cite[Theorem~11.3]{MR3452186}).
\end{Remark}

The notion of atoms in the category of modules over a ring was introduced by Storrer \cite{MR0360717} in order to reformulate the theory of primary decomposition due to Goldman \cite{MR0245608}. The definition below is a generalized version, which was given in \cite{MR2964615} to deal with arbitrary abelian categories.

\begin{Definition}\label{2309230982}
	Let $\cG$ be a Grothendieck category.
	\begin{enumerate}
		\item\label{2309209841} A nonzero object $H$ in $\cG$ is called \emph{monoform} if for every nonzero subobject $L$ of $H$, the only subobject of $H$ that is isomorphic to some subobject of $H/L$ is the zero subobject.
		\item\label{0120918234} We say that two monoform objects $H_{1}$ and $H_{2}$ are \emph{atom-equivalent} if $H_{1}$ has a nonzero subobject that is isomorphic to some subobject of $H_{2}$.
		\item\label{1092312033} The \emph{atom spectrum} of $\cG$ is defined to be
		\begin{equation*}
			\ASpec\cG:=\frac{\setwithtext{monoform objects in $\cG$}}{\text{atom equivalence}}.
		\end{equation*}
		The equivalence class of a monoform object $H$ is denoted by $\overline{H}$. Each element of $\ASpec\cG$ is called an \emph{atom} in $\cG$.
	\end{enumerate}
\end{Definition}

The atom spectrum has a structure of a topological space and a partially ordered set:

\begin{Definition}\label{0989183411}
	Let $\cG$ be a Grothendieck category.
	\begin{enumerate}
		\item\label{0912908238} For each object $M$ in $\cG$, we define the \emph{atom support} of $M$ to be
		\begin{equation*}
			\ASupp M:=\setwithcondition{\overline{H}\in\ASpec\cG}{\textnormal{$H$ is a monoform subquotient of $M$}}.
		\end{equation*}
		\item\label{2890231090} We define the \emph{localizing topology} on $\ASpec\cG$ as follows: A subset $\Phi$ of $\ASpec\cG$ is open if and only if $\Phi=\ASupp M$ for some $M\in\cG$.
		\item\label{0912902133} We define a partial order $\leq$ on $\ASpec\cG$, called the \emph{specialization order}, as follows: For $\alpha,\beta\in\ASpec\cG$, $\alpha\leq\beta$ if and only if $\alpha$ belongs to the closure of $\set{\beta}$ with respect to the localizing topology.
	\end{enumerate}
\end{Definition}

The localizing topology satisfies the axioms of a topology (see \cite[Proposition~3.2]{MR3351569}). The atom spectrum endowed with the localizing topology is a Kolmogorov space, and this implies that the binary relation $\leq$ satisfies the axioms of a partial order (see \cite[section~4]{MR3351569}).

\begin{Remark}\label{0913290123}
	For a commutative ring $R$, there is a bijection $\Spec R\isoto\ASpec(\Mod R)$ given by $\kp\mapsto\overline{R/\kp}$ as shown in \cite[p.~631]{MR0360717}. This is an isomorphism of partially ordered sets (\cite[Proposition~4.3]{MR3351569}). A subset of $\ASpec(\Mod R)$ is open with respect to the localizing topology if and only if the corresponding subset of $\Spec R$ is closed under specialization, that is, it is upward-closed with respect to the inclusion of prime ideals (\cite[Proposition~7.2 (2)]{MR2964615}). Thus the localizing topology is different from the Zariski topology.
	
	For a Grothendieck category $\cG$, every open subset of $\ASpec\cG$ is upward-closed with respect to the specialization order, but the converse does not hold in general (see \cite[Remark~2.11]{arXiv:1711.06946}).
\end{Remark}

\begin{Remark}\label{0920892345}
	For a locally noetherian Grothendieck category $\cG$, there is a canonical bijection between the atom spectrum of $\cG$ and the \emph{Gabriel spectrum} of $\cG$, which consists of all isomorphism classes of indecomposable injective objects in $\cG$. See \cite[Proposition~2.3]{arXiv:1711.06946} for the definitions of the bijection and the corresponding partial order on the Gabriel spectrum.
\end{Remark}

\begin{Remark}\label{0198241233}
	If $\cX$ is a weakly closed subcategory of a Grothendieck category $\cG$, then there is a canonical injective continuous map $\ASpec\cX\into\ASpec\cG$ defined by $\overline{H}\mapsto\overline{H}$. This map is also a homomorphism of partially ordered sets, and induces a homeomorphism from $\ASpec\cX$ to an open subset of $\ASpec\cG$, which was denoted by $\ASupp\cX$ in \cite[Proposition~5.12]{MR3351569}. We identify $\ASpec\cX$ with its image in $\ASpec\cG$.
\end{Remark}

\begin{Definition}\label{0923109234}
	Let $\cG$ be a Grothendieck category. An object $H$ in $\cG$ is called \emph{compressible} if each nonzero subobject of $H$ contains a subobject isomorphic to $H$.
\end{Definition}

\begin{Remark}\label{0981081238}
	In a locally noetherian Grothendieck category, every compressible object is monoform (see \cite[Proposition~2.12 (2)]{MR3922832}). On the other hand, a monoform object does not necessarily have a compressible object. Indeed, Goodearl \cite{MR588447} gave an example of a left and right noetherian ring $R$ that admits a monoform right $R$-module with no compressible submodule.
\end{Remark}

\begin{Remark}\label{0982234906}
	Let $\cG$ be a Grothendieck category with enough compressible objects. Then each atom $\alpha$ in $\cG$ is represented by some compressible object $H$, and $\ASupp H$ is the \emph{smallest} open subset of $\ASpec\cG$ among those containing $\alpha$ (see \cite[Definition~3.1 and Proposition~3.2]{MR3351569}). This implies that the intersection of an arbitrary family of open subsets of $\ASpec\cG$ is again an open subset. Such a topological space is called an \emph{Alexandroff space}.
	
	As a consequence, the topology of $\ASpec\cG$ can be recovered from its partial order as follows: a subset of $\ASpec\cG$ is open if and only if it is upward-closed with respect to the partial order (see, for example, \cite[Proposition~4.1]{MR3351569}). In particular, the Grothendieck categories that we will construct in \cref{0893340933} have this property.
\end{Remark}

\section{Grothendieck categories associated to colored quivers}
\label{0295378257}

We fix a field $K$. For a given colored quiver $\Gamma$, we will associate a $K$-algebra $F_{\Gamma}$, an $F_{\Gamma}$-module $M_{\Gamma}$, and a Grothendieck category $\cG_{\Gamma}$.

For a set $C$, we denote by $C^{*}$ the free (multiplicative) monoid on $C$. The identity element is denoted by $1$, and any other element of $C^{*}$ is of the form $\bc=c_{1}c_{2}\cdots c_{l}$, where $l\geq 1$ and $c_{1},\ldots,c_{l}\in C$.

\begin{Definition}\label{0445398352}\leavevmode
	\begin{enumerate}
		\item\label{0278952473} A \emph{colored quiver} is a sextuple $\Gamma=(\Gamma_{0},\Gamma_{1},C,s,t,u)$, where $\Gamma_{0}$, $\Gamma_{1}$, and $C$ are sets, and
		\begin{equation*}
			s\colon\Gamma_{1}\to\Gamma_{0},\quad t\colon\Gamma_{1}\to\Gamma_{0},\quad\text{and}\quad u\colon\Gamma_{1}\to C
		\end{equation*}
		are maps. The elements of $\Gamma_{0}$, $\Gamma_{1}$, and $C$ are called \emph{vertices}, \emph{arrows}, and \emph{colors} in $\Gamma$, respectively. For each arrow $r$ in $\Gamma$, the images $s(r)$, $t(r)$, and $u(r)$ are called the \emph{source}, \emph{target}, and \emph{color} of $r$, respectively.
		
		We often write $\Gamma=(\Gamma_{0},\Gamma_{1},C)$ by omitting $s$, $t$, and $u$.
		\item\label{0981403253} Let $\Gamma=(\Gamma_{0},\Gamma_{1},C)$ be a colored quiver. A \emph{path} of length $l\geq 1$ in $\Gamma$ is a sequence $\br=r_{1}r_{2}\cdots r_{l}$ of $l$ arrows in $\Gamma$ such that $t(r_{i})=s(r_{i+1})$ for all $i=1,\ldots,l-1$. A \emph{path} of length $0$ is a symbol $e_{v}$, where $v$ is a vertex in $\Gamma$.
		
		For each path $\br=r_{1}r_{2}\cdots r_{l}$ in $\Gamma$, the \emph{source} and \emph{target} of $\br$ are defined to be
		\begin{equation*}
			s(\br):=s(r_{1})\quad\text{and}\quad t(\br):=t(r_{l}),
		\end{equation*}
		respectively. The \emph{sequence of colors} of $\br$ is
		\begin{equation*}
			u(\br):=u(r_{1})u(r_{2})\cdots u(r_{l})\in C^{*}.
		\end{equation*}
		For $\br=e_{v}$, let $s(\br)=t(\br)=v$ and $u(\br)=1\in C^{*}$.
	\end{enumerate}
\end{Definition}

\begin{Notation}\label{0184739012}
	An arrow $r$ in a colored quiver is visualized as
	\begin{equation*}
		\begin{tikzcd}
			s(r)\ar[r,"u(r)"] & t(r)
		\end{tikzcd}.
	\end{equation*}
	Thus, for example, a path $\br=r_{1}r_{2}$ of length $2$ can be written as
	\begin{equation*}
		\begin{tikzcd}
			s(r_{1})\ar[r,"u(r_{1})"] & v\ar[r,"u(r_{2})"] & t(r_{2})
		\end{tikzcd},
	\end{equation*}
	where $v:=t(r_{1})=s(r_{2})$.
\end{Notation}

\begin{Definition}\label{2358974250}
	Let $C$ be a set.
	\begin{enumerate}
		\item\label{0798423234} We define the $K$-algebra $\KC$ as follows: As a $K$-vector space,
		\begin{equation*}
			\KC=\prod_{\bc\in C^{*}}K\bc,
		\end{equation*}
		and each element is written as a formal sum $\sum_{\bc\in C^{*}}\lambda_{\bc}\bc$ with $\lambda_{\bc}\in K$. The multiplication is defined as
		\begin{equation*}
			\Bigg(\sum_{\bc\in C^{*}}\lambda_{\bc}\bc\Bigg)\Bigg(\sum_{\bc\in C^{*}}\lambda'_{\bc}\bc\Bigg)=\sum_{\bc\in C^{*}}\Bigg(\sum_{\bc_{1},\bc_{2}}\lambda_{\bc_{1}}\lambda'_{\bc_{2}}\Bigg)\bc,
		\end{equation*}
		where $\bc_{1}$ and $\bc_{2}$ run over all pairs of elements of $C^{*}$ with $\bc_{1}\bc_{2}=\bc$. Note that there are only finitely many such pairs for each $\bc$.
		\item\label{3203823847} For each $f=\sum_{\bc\in C^{*}}\lambda_{\bc}\bc\in\KC$, define its \emph{support} to be
		\begin{equation*}
			\supp f:=\setwithcondition{\bc\in C^{*}}{\lambda_{\bc}\neq 0}.
		\end{equation*}
	\end{enumerate}
\end{Definition}

\begin{Definition}\label{0294572409}
	Let $\Gamma=(\Gamma_{0},\Gamma_{1},C)$ be a colored quiver.
	\begin{enumerate}
		\item\label{0295284344} For each subset $B\subset C^{*}$ and each vertex $v$ in $\Gamma$, define the set $P_{v}(B)$ to be the set of all paths $\br$ in $\Gamma$ such that $t(\br)=v$ and $u(\br)\in B$.
		\item\label{2345987230} We say that an element $f\in\KC$ is \emph{admissible} if $P_{v}(\supp f)$ is a finite set for all vertices $v$ in $\Gamma$.
		\item\label{0923485423} Denote by $F_{\Gamma}$ the set of all admissible elements of $\KC$.
	\end{enumerate}
\end{Definition}

\begin{Lemma}\label{0897423522}
	Let $\Gamma=(\Gamma_{0},\Gamma_{1},C)$ be a colored quiver. Then $F_{\Gamma}$ is a $K$-subalgebra of $\KC$.
\end{Lemma}

\begin{proof}
	For each $B$ and $v$, we have
	\begin{equation*}
		P_{v}(B)=\bigcup_{\bc\in B}P_{v}(\set{\bc}).
	\end{equation*}
	Thus the inclusions
	\begin{equation*}
		\supp(\lambda f)\subset\supp f\quad and\quad\supp(f+g)\subset\supp f\cup\supp g,
	\end{equation*}
	where $f,g\in\KC$ and $\lambda\in K$, imply that $F_{\Gamma}$ is a $K$-subspace of $\KC$. In order to conclude the claim, it suffices to show that $fg$ is admissible whenever $f,g\in\KC$ are admissible.
	
	Let $v$ be a vertex in $\Gamma$. Since $P_{v}(\supp g)$ is finite, we can write $P_{v}(\supp g)=\set{\br_{1},\ldots,\br_{n}}$. Take arbitrary $\br\in P_{v}(\supp fg)$. Since $u(\br)\in\supp fg$, there exist paths $\br'$ and $\br''$ in $\Gamma$ such that $\br=\br'\br''$, $u(\br')\in\supp f$, and $u(\br'')\in\supp g$. Here $\br''\in P_{v}(\supp g)$, so that $\br''=\br_{i}$ for some $i$, and $\br'\in P_{s(\br_{i})}(\supp f)$. We have shown that every path in $P_{v}(\supp fg)$ is of the form $\br'\br''$, where $\br'$ belongs to the finite set $P_{s(\br_{i})}(\supp f)$ for some $i=1,\ldots,n$ and $\br''$ belongs to the finite set $P_{v}(\supp g)$. Therefore $P_{v}(\supp fg)$ is a finite set. $fg$ is admissible.
\end{proof}

\begin{Definition}\label{9720252347}
	Let $\Gamma=(\Gamma_{0},\Gamma_{1},C)$ be a colored quiver. We define the right $F_{\Gamma}$-module $M_{\Gamma}$ as follows: As a $K$-vector space,
	\begin{equation*}
		M_{\Gamma}=\prod_{v\in\Gamma_{0}}Kx_{v},
	\end{equation*}
	where $Kx_{v}$ is a one-dimensional space generated by $x_{v}$. Each element of $M_{\Gamma}$ is written as a formal sum $\sum_{v\in\Gamma_{0}}\mu_{v}x_{v}$. For each
	\begin{equation*}
		f=\sum_{\bc\in C^{*}}\lambda_{\bc}\bc\in F_{\Gamma}\quad\text{and}\quad y=\sum_{v\in\Gamma_{0}}\mu_{v}x_{v}\in M_{\Gamma},
	\end{equation*}
	define
	\begin{equation*}
		yf:=\sum_{\br}\lambda_{u(\br)}\mu_{s(\br)}x_{t(\br)},
	\end{equation*}
	where $\br$ runs over all paths in $\Gamma$. The sum makes sense because, for each vertex $v$ in $\Gamma$, there are only finitely many paths $\br$ satisfying $t(\br)=v$ and $\lambda_{u(\br)}\neq 0$ by the admissibility of $f$. It is easy to verify that $M_{\Gamma}$ is actually a right $F_{\Gamma}$-module.
	
	For each $y=\sum_{v\in\Gamma_{0}}\mu_{v}x_{v}\in M_{\Gamma}$, define its \emph{support} to be
	\begin{equation*}
		\supp y:=\setwithcondition{v\in\Gamma_{0}}{\mu_{v}\neq 0}.
	\end{equation*}
\end{Definition}

\begin{Definition}\label{0298337452}
	Let $\Gamma=(\Gamma_{0},\Gamma_{1},C)$ be a (set-theoretically small) colored quiver. We define $\cG_{\Gamma}$ to be the smallest weakly closed subcategory of $\Mod F_{\Gamma}$ containing $M_{\Gamma}$.
\end{Definition}

\begin{Example}\label{1032489742}
	If $\Gamma=(\Gamma_{0},\Gamma_{1},C)$ is the colored quiver
	\begin{equation*}
		\begin{tikzcd}
			v\ar[out=240,in=300,loop,looseness=5,"c"']
		\end{tikzcd}
	\end{equation*}
	with $C=\set{c}$, then $\KC$ is the (commutative) formal power series ring $K\llbracket c\rrbracket$ and $F_{\Gamma}$ is the polynomial ring $K[c]$. The $K[c]$-module $M_{\Gamma}$ is isomorphic to the simple module $K[c]/(c-1)K[c]$. Thus $\cG_{\Gamma}=\Mod (K[c]/(c-1)K[c])\cong\Mod K$.
\end{Example}

\section{Construction of Grothendieck categories}
\label{0983897234}

As in the previous section, let $K$ be a field. We denote by $\bbN$ the set of nonnegative integers.

We will introduce an operations to obtain a new colored quiver from a given one, which allows us to construct a colored quiver $\Gamma$ such that the Grothendieck category $\cG_{\Gamma}$ satisfies the desired properties in \cref{0984767827}.

\begin{Definition}\label{8790234583}
	Let $\Gamma=(\Gamma_{0},\Gamma_{1},C)$ be a colored quiver. Define a new colored quiver $\widetilde{\Gamma}=(\widetilde{\Gamma}_{0},\widetilde{\Gamma}_{1},\widetilde{C})$ as follows:
	\begin{itemize}
		\item $\widetilde{\Gamma}_{0}=\bbN\times\Gamma_{0}$.
		\item $\widetilde{\Gamma}_{1}=(\bbN\times\Gamma_{1})\amalg\setwithcondition{r^{i}_{v,w}}{\textnormal{$i\in\bbN$, $v,w\in\Gamma_{0}$}}$, where $\amalg$ denotes the disjoint union of sets and $r^{i}_{v,w}$ are pairwise distinct symbols.
		\item $\widetilde{C}=C\amalg\setwithcondition{c_{v,w}}{v,w\in\Gamma_{0}}$, where $c_{v,w}$ are pairwise distinct symbols.
		\item For each $(i,r)\in\bbN\times\Gamma_{1}$,
		\begin{equation*}
			s(i,r)=(i,s(r)),\quad t(i,r)=(i,t(r)),\quad\text{and}\quad u(i,r)=u(r).
		\end{equation*}
		\item For each $r^{i}_{v,w}$,
		\begin{equation*}
			s(r^{i}_{v,w})=(i,v),\quad t(r^{i}_{v,w})=(i+1,w),\quad\text{and}\quad u(r^{i}_{v,w})=c_{v,w}.
		\end{equation*}
	\end{itemize}
\end{Definition}

\begin{Remark}\label{0910983243}
	\cref{8790234583} is actually a special case of \cite[Definition~7.17]{MR3351569}, but it is different from the operation observed in \cite[Proposition~7.22]{MR3351569}, which played a crucial role in the proof of \cite[Theorem~7.23]{MR3351569}.
	
	Indeed, if we set $\Omega$ to be the colored quiver
	\begin{equation*}
		\begin{tikzcd}
			\omega_{0}\ar[r,"\xi"] & \omega_{1}\ar[r,"\xi"] & \cdots
		\end{tikzcd}
	\end{equation*}
	with the set of colors $\Xi=\set{\xi}$ and all $\Gamma^{\omega_{i}}$ to be $\Gamma$ in \cite[Definition~7.17]{MR3351569}, then the output $\Omega(\set{\Gamma^{\omega}}_{\omega\in\Omega_{0}})$ is the $\widetilde{\Gamma}$ in \cref{8790234583} (up to change of symbols).
	
	On the other hand, in order to obtain the colored quiver $\Gamma$ in \cite[Proposition~7.22]{MR3351569}, we have to set $\Omega$ to be
	\begin{equation*}
		\begin{tikzcd}
			\omega_{0}\ar[r,"\xi_{0}"] & \omega_{1}\ar[r,"\xi_{1}"] & \cdots
		\end{tikzcd}
	\end{equation*}
	with the set of colors $\Xi=\setwithcondition{\xi_{i}}{i\in\bbN}$, where $\xi_{i}$ are pairwise distinct, and each $\Gamma^{\omega_{i}}$ to be $\Gamma^{i}$.
\end{Remark}

\begin{Example}\label{9287345224}
	If $\Gamma=(\Gamma_{0},\Gamma_{1},C)$ is the colored quiver
	\begin{equation*}
		\begin{tikzcd}
			v\ar[d,"c"] \\
			w
		\end{tikzcd}
	\end{equation*}
	with $C=\set{c}$, then $\widetilde{\Gamma}$ is
	\begin{equation*}
		\begin{tikzcd}[row sep=15mm, column sep=15mm]
			(0,v)\ar[d,"c"]\ar[r,"c_{v,v}"]\ar[dr,"c_{v,w}" description, pos=0.3] & (1,v)\ar[d,"c"]\ar[r,"c_{v,v}"]\ar[dr,"c_{v,w}" description, pos=0.3] & (2,v)\ar[d,"c"]\ar[r,"c_{v,v}"]\ar[dr,"c_{v,w}" description, pos=0.3] & \cdots \\
			(0,w)\ar[ur,"c_{w,v}" description, pos=0.3]\ar[r,"c_{w,w}"'] & (1,w)\ar[ur,"c_{w,v}" description, pos=0.3]\ar[r,"c_{w,w}"'] & (2,w)\ar[ur,"c_{w,v}" description, pos=0.3]\ar[r,"c_{w,w}"'] & \cdots
		\end{tikzcd}
	\end{equation*}
	with $\widetilde{C}=\set{c,\,c_{v,v},\,c_{v,w},\,c_{w,v},\,c_{w,w}}$.
\end{Example}

\begin{Example}\label{0894320823}
	If $\Gamma$ is the colored quiver with only one vertex $v$, no arrows, and no colors, then $\widetilde{\Gamma}$ is
	\begin{equation*}
		\begin{tikzcd}
			(0,v)\ar[r,"c_{v,v}"] & (1,v)\ar[r,"c_{v,v}"] & (2,v)\ar[r,"c_{v,v}"] & \cdots
		\end{tikzcd}
	\end{equation*}
	with $\widetilde{C}=\set{c_{v,v}}$. Thus $\KtC$ is the formal power series ring $K\llbracket c_{v,v}\rrbracket$ with a single variable $c_{v,v}$, and $M_{\Gamma}=F_{\Gamma}=\KtC=K\llbracket c_{v,v}\rrbracket$. Therefore $\cG_{\Gamma}=\Mod K\llbracket c_{v,v}\rrbracket$.
\end{Example}

\begin{Remark}\label{0892340924}
	In the setting of \cref{8790234583}, $\KC$ is a $K$-subalgebra of $\KtC$ since $C\subset\widetilde{C}$. Write the inclusion map as $\nu\colon\KC\into\KtC$. We can define a surjective ring homomorphism $\pi\colon\KtC\onto\KC$ by
	\begin{equation*}
		\sum_{\bc\in\widetilde{C}^{*}}\lambda_{\bc}\bc\mapsto\sum_{\bc\in C^{*}}\lambda_{\bc}\bc,
	\end{equation*}
	which satisfies $\pi\circ\nu=\id$. It can be verified that $\nu$ and $\pi$ induce $K$-algebra homomorphisms between $F_{\Gamma}$ and $F_{\widetilde{\Gamma}}$. The kernel of the induced surjective homomorphism $F_{\widetilde{\Gamma}}\onto F_{\Gamma}$ is
	\begin{equation*}
		I_{\Gamma}:=\setwithcondition[\bigg]{\sum_{\bc\in\widetilde{C}^{*}}\lambda_{\bc}\bc\in F_{\widetilde{\Gamma}}}{\textnormal{$\lambda_{\bc}=0$ for all $\bc\in C^{*}$}}.
	\end{equation*}
	Therefore
	\begin{equation*}
		\Mod F_{\Gamma}\isoto\Mod\frac{F_{\widetilde{\Gamma}}}{I_{\Gamma}}=\setwithcondition{N\in\Mod F_{\widetilde{\Gamma}}}{NI_{\Gamma}=0}.
	\end{equation*}
	We regard $\Mod F_{\Gamma}$ as a closed subcategory of $\Mod F_{\widetilde{\Gamma}}$ in this way.
\end{Remark}

\begin{Definition}\label{0923489043}
	Let $\Gamma=(\Gamma_{0},\Gamma_{1},C)$ be a colored quiver and let $\widetilde{\Gamma}=(\widetilde{\Gamma}_{0},\widetilde{\Gamma}_{1},\widetilde{C})$ as in \cref{8790234583}. For each $i\in\bbN$, define the $F_{\widetilde{\Gamma}}$-submodule $M_{\geq i}\subset M_{\widetilde{\Gamma}}$ by
	\begin{equation*}
		M_{\geq i}=\setwithcondition[\bigg]{y\in M_{\widetilde{\Gamma}}}{\supp y\subset\bigcup_{j\geq i}(\set{j}\times\Gamma_{0})}.
	\end{equation*}
	Define the $K$-subspace $M_{i}\subset M_{\widetilde{\Gamma}}$ by
	\begin{equation*}
		M_{i}=\setwithcondition{y\in M_{\widetilde{\Gamma}}}{\supp y\subset(\set{i}\times\Gamma_{0})}.
	\end{equation*}
\end{Definition}

\begin{Remark}\label{0919081802}
	In the setting of \cref{0923489043}, there is an isomorphism $M_{\widetilde{\Gamma}}\isoto M_{\geq i}$ of $F_{\widetilde{\Gamma}}$-modules given by
	\begin{equation*}
		\sum_{j\in\bbN}\sum_{v\in\Gamma_{0}}\lambda_{(j,v)}x_{(j,v)}\mapsto\sum_{j\in\bbN}\sum_{v\in\Gamma_{0}}\lambda_{(j,v)}x_{(j+i,v)}.
	\end{equation*}
	
	We identify the $K$-vector space $M_{i}$ with the right $F_{\widetilde{\Gamma}}$-module $M_{\geq i}/M_{\geq i+1}$ via the composite
	\begin{equation*}
		M_{i}\into M_{\geq i}\onto M_{\geq i}/M_{\geq i+1}.
	\end{equation*}
	Note that $M_{i}=M_{\geq i}/M_{\geq i+1}$ belongs to the closed subcategory $\Mod F_{\Gamma}$ of $\Mod F_{\widetilde{\Gamma}}$. We can define an isomorphism $M_{\Gamma}\isoto M_{i}$ of $F_{\widetilde{\Gamma}}$-modules (also of $F_{\Gamma}$-modules) by
	\begin{equation*}
		\sum_{v\in\Gamma_{0}}\lambda_{v}x_{v}\mapsto\sum_{v\in\Gamma_{0}}\lambda_{v}x_{(i,v)}.
	\end{equation*}
	This means that $M_{\Gamma}$ is isomorphic to a subquotient of $M_{\widetilde{\Gamma}}$ in $\Mod F_{\widetilde{\Gamma}}$. Consequently, $\cG_{\Gamma}$ is a weakly closed subcategory of $\cG_{\widetilde{\Gamma}}$.
\end{Remark}

The following lemmas will be used to show some properties of the Grothendieck category $\cG_{\widetilde{\Gamma}}$:

\begin{Lemma}\label{0987342584}
	Let $\Gamma=(\Gamma_{0},\Gamma_{1},C)$ be a colored quiver and let $\widetilde{\Gamma}=(\widetilde{\Gamma}_{0},\widetilde{\Gamma}_{1},\widetilde{C})$ as in \cref{8790234583}. Let $y\in M_{\widetilde{\Gamma}}$, and take $i\in\bbN$ such that
	\begin{equation*}
		(\supp y)\cap(\set{i}\times\Gamma_{0})\neq\emptyset.
	\end{equation*}
	Then $M_{\geq i+1}\subset yF_{\widetilde{\Gamma}}$.
\end{Lemma}

\begin{proof}
	We can assume $y\in M_{\geq i}$ by replacing $i$ by the smallest one without loss of generality. Write
	\begin{equation*}
		y=\sum_{j\geq i}\sum_{v\in\Gamma_{0}}\lambda_{(j,v)}x_{(j,v)}.
	\end{equation*}
	By the assumption, there exists $u\in\Gamma_{0}$ such that $\lambda_{(i,u)}\neq 0$. By replacing $y$ by $\lambda_{(i,u)}^{-1}y$, we may assume $\lambda_{(i,u)}=1$.
	
	Let $z\in M_{\geq i+1}$. For each integer $d\geq 1$, define $z_{d}\in M_{\geq i+d}$ and $f_{d}\in F_{\widetilde{\Gamma}}$ inductively as follows: Let $z_{1}:=z$. When $z_{d}\in M_{\geq i+d}$ is defined, write
	\begin{equation*}
		z_{d}=\sum_{j\geq i+d}\sum_{v\in\Gamma_{0}}\mu_{(j,v)}x_{(j,v)}.
	\end{equation*}
	and define
	\begin{equation*}
		f_{d}:=\sum_{v\in\Gamma_{0}}\mu_{(i+d,v)}c_{(u,u)}^{d-1}c_{(u,v)}.
	\end{equation*}
	Let $z_{d+1}:=z_{d}-yf_{d}$. Since
	\begin{equation*}
		\begin{split}
			yf_{d}
			&=\sum_{j\geq i}\sum_{v\in\Gamma_{0}}\lambda_{(j,u)}\mu_{(i+d,v)}x_{(j+d,v)}
			=\sum_{j\geq i+d}\sum_{v\in\Gamma_{0}}\lambda_{(j-d,u)}\mu_{(i+d,v)}x_{(j,v)}\\
			&=\sum_{v\in\Gamma_{0}}\mu_{(i+d,v)}x_{(i+d,v)}+\sum_{j\geq i+d+1}\sum_{v\in\Gamma_{0}}\lambda_{(j-d,u)}\mu_{(i+d,v)}x_{(j,v)},
		\end{split}
	\end{equation*}
	we have $z_{d+1}\in M_{\geq i+d+1}$. It follows that $f:=\sum_{d\geq 1}f_{d}\in F_{\widetilde{\Gamma}}$, and
	\begin{equation*}
		z-yf=z_{1}-\sum_{d\geq 1}yf_{d}=0.
	\end{equation*}
	Thus $z=yf\in yF_{\widetilde{\Gamma}}$.
\end{proof}

\begin{Lemma}\label{0897328934}
	Let $\Gamma=(\Gamma_{0},\Gamma_{1},C)$ be a colored quiver and let $\widetilde{\Gamma}=(\widetilde{\Gamma}_{0},\widetilde{\Gamma}_{1},\widetilde{C})$ as in \cref{8790234583}. For each nonzero $F_{\widetilde{\Gamma}}$-submodule $L\subset M_{\widetilde{\Gamma}}$, there exists $i\in\bbN$ such that
	\begin{equation*}
		L=(L\cap M_{i})\oplus M_{\geq i+1}
	\end{equation*}
	as a $K$-vector space.
\end{Lemma}

\begin{proof}
	Let $i$ be the smallest number satisfying
	\begin{equation*}
		(\supp y)\cap(\set{i}\times\Gamma_{0})\neq\emptyset
	\end{equation*}
	for some $y\in L$. Then by \cref{0987342584},
	\begin{equation*}
		M_{\geq i+1}\subset yF_{\widetilde{\Gamma}}\subset L\subset M_{\geq i}.
	\end{equation*}
	In particular, $(L\cap M_{i})\oplus M_{\geq i+1}\subset L$.
	
	Each $z\in L$ can be written as $z=z'+z''$, where $z'\in M_{i}$ and $z''\in M_{\geq i+1}$. There exists $f\in F_{\widetilde{\Gamma}}$ such that $yf=z''$. Hence $z'=z-yf\in L\cap M_{i}$. This shows that $L\subset (L\cap M_{i})\oplus M_{\geq i+1}$.
\end{proof}

\begin{Proposition}\label{3409359033}
	Let $\Gamma=(\Gamma_{0},\Gamma_{1},C)$ be a colored quiver with $\Gamma_{0}\neq\emptyset$ and let $\widetilde{\Gamma}=(\widetilde{\Gamma}_{0},\widetilde{\Gamma}_{1},\widetilde{C})$ as in \cref{8790234583}. Suppose that $\Gamma$ satisfies the following conditions:
	\begin{enumerate}\renewcommand{\theenumi}{\alph{enumi}}
		\item\label{0920934343} $M_{\Gamma}$ is a noetherian object in $\cG_{\Gamma}$.
		\item\label{2390823091} $\cG_{\Gamma}$ has enough compressible objects, that is, every nonzero object in $\cG_{\Gamma}$ has a compressible subobject.
	\end{enumerate}
	Then the following hold:
	\begin{enumerate}
		\item\label{8093420983} $M_{\widetilde{\Gamma}}$ is a compressible noetherian object in $\cG_{\widetilde{\Gamma}}$. Consequently, $\cG_{\widetilde{\Gamma}}$ is a locally noetherian Grothendieck category.
		\item\label{1089012347} $\cG_{\widetilde{\Gamma}}$ has enough compressible objects.
		\item\label{0329023890} As a partially ordered set,
		\begin{equation*}
			\ASpec\cG_{\widetilde{\Gamma}}=\ASpec\cG_{\Gamma}\cup\set[\Big]{\overline{M_{\widetilde{\Gamma}}}},
		\end{equation*}
		where $\overline{M_{\widetilde{\Gamma}}}$ is smaller than any element in $\ASpec\cG_{\Gamma}$.
	\end{enumerate}
\end{Proposition}

\begin{proof}
	\cref{8093420983} Every nonzero submodule $L$ of $M_{\widetilde{\Gamma}}$ contains some $M_{\geq i+1}$ as a subobject by \cref{0897328934}. Since $M_{\geq i+1}$ is isomorphic to $M_{\widetilde{\Gamma}}$ by \cref{0919081802}, $M_{\widetilde{\Gamma}}$ is compressible. Since we have the filtration
	\begin{equation*}
		M_{\widetilde{\Gamma}}=M_{\geq 0}\supset M_{\geq 1}\supset\cdots,
	\end{equation*}
	it suffices to show that $M_{\geq i}/M_{\geq i+1}$ is a noetherian for each $i\geq 0$. Again by \cref{0919081802}, this is isomorphic to $M_{\Gamma}$. Since $M_{\Gamma}$ is a noetherian $F_{\Gamma}$-module by the assumption, it is also noetherian as a $M_{\widetilde{\Gamma}}$-module.
	
	As explained in the paragraph before Remark~7.9 in \cite{MR3351569}, it follows that $\cG_{\widetilde{\Gamma}}$ is locally noetherian.
	
	\cref{1089012347} Since $\cG_{\widetilde{\Gamma}}$ is locally noetherian, every nonzero object contains a monoform subobject (\cite[Theorem~2.9]{MR2964615}). As in the paragraph after Remark~7.8 in \cite{MR3351569}, $\ASupp M_{\widetilde{\Gamma}}=\ASpec\cG_{\widetilde{\Gamma}}$. Thus each monoform object contains a subobject that is isomorphic to a monoform object of the form $L'/L$, where $L\subset L'$ are subobjects of $M_{\widetilde{\Gamma}}$.
	
	If $L=0$, then $L'$ contains some $M_{\geq i+1}$ as a subobject, which is compressible since $M_{\geq i+1}\cong M_{\widetilde{\Gamma}}$.
	
	If $L\neq 0$, then $L=(L\cap M_{i})\oplus M_{\geq i+1}$ and $L\cap M_{i}\subsetneq M_{i}$ for some $i$, and hence $L'/L$ contains the submodule
	\begin{equation*}
		\frac{(L'\cap M_{i})\oplus M_{\geq i+1}}{(L\cap M_{i})\oplus M_{\geq i+1}}\cong\frac{L'\cap M_{i}}{L\cap M_{i}},
	\end{equation*}
	which is a nonzero subquotient of $M_{i}\cong M_{\Gamma}\in\cG_{\Gamma}$. Hence it contains a compressible subobject by the assumption.
	
	\cref{0329023890} The proof of \cref{1089012347} shows in particular that every monoform object in $\cG_{\widetilde{\Gamma}}$ contains either a submodule isomorphic to $M_{\widetilde{\Gamma}}$ or a compressible (thus monoform) subobject in $\cG_{\Gamma}$. Hence the equality in the claim follows. Since $M_{\widetilde{\Gamma}}$ is a compressible object and its atom support is $\ASpec\cG_{\widetilde{\Gamma}}$, the atom $\overline{M_{\widetilde{\Gamma}}}$ is smallest by \cite[Proposition~4.2]{MR3351569}. Since no nonzero subobject of $M_{\widetilde{\Gamma}}$ belongs to $\cG_{\Gamma}$, the atom $\overline{M_{\widetilde{\Gamma}}}$ does not belong to $\ASpec\cG_{\Gamma}$.
\end{proof}

We also use the disjoint union of colored quivers:

\begin{Definition}[{\cite[Definition~7.17]{MR3351569}}]\label{3429083203}
	Let $\set{\Gamma^{i}}_{i\in I}$ be a family of colored quivers, where $\Gamma^{i}=(\Gamma^{i}_{0},\Gamma^{i}_{1},C^{i})$. Its \emph{disjoint union} is defined to be
	\begin{equation*}
		\coprod_{i\in I}\Gamma^{i}=\bigg(\coprod_{i\in I}\Gamma^{i}_{0},\ \coprod_{i\in I}\Gamma^{i}_{1},\ \bigcup_{i\in I}C^{i}\bigg),
	\end{equation*}
	where $s$, $t$, and $u$ are defined to be those induced from the colored quivers $\Gamma^{i}$.
\end{Definition}

\begin{Remark}\label{0923459823}
	In the setting of \cref{3429083203}, each $\Mod\Gamma^{i}$ is regarded as a closed subcategory of $\Mod\Gamma$, where $\Gamma:=\coprod_{i\in I}\Gamma^{i}$, in the same way as \cref{0892340924}. Thus each $\cG_{\Gamma^{i}}$ is a weakly closed subcategory of $\cG_{\Gamma}$.
\end{Remark}

\begin{Proposition}\label{0892392343}
	Let $\set{\Gamma^{i}}_{i=1}^{n}$ be a finite family of colored quivers and let $\Gamma:=\coprod_{i=1}^{n}\Gamma^{i}$. Suppose that each $\Gamma_{i}$ satisfies the following conditions:
	\begin{enumerate}\renewcommand{\theenumi}{\alph{enumi}}
		\item\label{0982309234} $M_{\Gamma^{i}}$ is a noetherian object in $\cG_{\Gamma^{i}}$.
		\item\label{1908278454} $\cG_{\Gamma^{i}}$ has enough compressible objects.
	\end{enumerate}
	Then the following hold:
	\begin{enumerate}
		\item\label{9801278123} $M_{\Gamma}$ is a noetherian object in $\cG_{\Gamma}$. Consequently, $\cG_{\Gamma}$ is a locally noetherian Grothendieck category.
		\item\label{9014324121} $\cG_{\Gamma}$ has enough compressible objects.
		\item\label{0912380491} As a partially ordered set,
		\begin{equation*}
			\ASpec\cG_{\Gamma}=\bigcup_{i=1}^{n}\ASpec\cG_{\Gamma^{i}}.
		\end{equation*}
	\end{enumerate}
\end{Proposition}

\begin{proof}
	These can be shown similarly to \cref{3409359033}, using $M_{\Gamma}=M_{\Gamma^{1}}\oplus\cdots\oplus M_{\Gamma^{n}}$.
\end{proof}

We are ready to prove the main result of this paper.

\begin{Theorem}\label{0893340933}
	Let $P$ be a finite partially ordered set. Then there exists a colored quiver $\Gamma$ satisfying the following properties:
	\begin{enumerate}
		\item\label{0934202345} $M_{\Gamma}$ is a noetherian object in $\cG_{\Gamma}$. Consequently, $\cG_{\Gamma}$ is a locally noetherian Grothendieck category.
		\item\label{9287409573} $\cG$ has enough compressible objects.
		\item\label{0290389513} $\ASpec\cG_{\Gamma}$ is isomorphic to $P$ as a partially ordered set.
	\end{enumerate}
\end{Theorem}

\begin{proof}
	For each $p\in P$, we associate a colored quiver $\Gamma^{p}$ inductively as follows:
	\begin{itemize}
		\item If $p$ is a maximal element, then $\Gamma^{p}=(\Gamma^{p}_{0},\Gamma^{p}_{1},C^{p})$ is
		\begin{equation*}
			\begin{tikzcd}
				v_{p}\ar[out=240,in=300,loop,looseness=5,"c_{p}"']
			\end{tikzcd}
		\end{equation*}
		with $C^{p}=\set{c_{p}}$, where $c_{p}$ is a symbol that is not introduced in the definition of any other $\Gamma^{q}$.
		\item If $p$ is not a maximal element, then let
		\begin{equation*}
			\Omega^{p}:=\coprod_{\substack{q\in P\\p<q}}\Gamma^{q}
		\end{equation*}
		and define $\Gamma^{p}:=\widetilde{\Omega^{p}}$ using \cref{8790234583}, where the symbols $c_{v,w}$ used in the definition of $\widetilde{\Omega^{p}}$ are chosen so that they are not introduced in the definition of any other $\Gamma^{q}$.
	\end{itemize}
	The colored quiver $\Gamma$ is defined by
	\begin{equation*}
		\Gamma:=\coprod_{p\in P}\Gamma^{p}.
	\end{equation*}
	Then all claims follow from \cref{1032489742}, \cref{3409359033}, and \cref{0892392343}.
\end{proof}



\end{document}